\documentclass{article}

\usepackage[utf8]{inputenc}
\usepackage{spconf}
\usepackage{mathtools}
\usepackage{tcolorbox}
\usepackage{amsmath}
\usepackage{amssymb}
\usepackage{mathtools}
\usepackage{algorithmic}
\usepackage{algorithm}
\usepackage{times}
\usepackage{tikz}
\usepackage{enumitem}

\input{mysymbol.sty}

\newcommand{\vect}{\mathrm{vec}}

\newtheorem{definition}{\hspace{0pt}\bf Definition}

\title{A Multi-resolution Low-rank Tensor Decomposition}

\name{\vspace{-.25cm}Sergio Rozada and Antonio G. Marques \thanks{Work supported by the Spanish NSF Grant SPGraph (PID2019-105032GB-I00), and by the Grants F661-MAPPING-UCI, F663-AAGNCS and F730-DSSP funded
		by the Comunidad de Madrid (CAM) and King Juan Carlos University
		(URJC).} 
}
\address{Dept. of Signal Theory and Communications, King Juan Carlos University, Madrid, Spain}

\begin{document}
\maketitle

\ninept

\begin{abstract}%
The (efficient and parsimonious) decomposition of higher-order tensors is a fundamental problem with numerous applications in a variety of fields. Several methods have been proposed in the literature to that end, with the Tucker and PARAFAC decompositions being the most prominent ones. Inspired by the latter, in this work we propose a multi-resolution low-rank tensor decomposition to describe (approximate) a tensor in a hierarchical fashion. The central idea of the decomposition is to recast the tensor into \emph{multiple} lower-dimensional tensors to exploit the structure at different levels of resolution. The method is first explained, an alternating least squares algorithm is discussed, and preliminary simulations illustrating the potential practical relevance are provided. 
\end{abstract}



\section{Introduction} \label{S:Introduction}

We live in a digital age where common-life devices, from smartphones to cars, generate massive amounts of data that provide researchers and practitioners a range of opportunities. Processing contemporary information comes, however, at a cost, since data sources are messy and heterogeneous. In this context, parsimonious models emerge as an ideal tool to enhance efficiency when processing such vast amounts of information. This can be done by leveraging the structure of the data, as is the case of information living in multiple (possibly many) dimensions. Multi-dimensional data are prevalent in numerous fields, with representative examples including chemometrics, bioengineering, communications, hyper-spectral imaging, or psychometrics \cite{bro1997parafac,cattell1944parallel}. Traditionally, matrices were used to model those datasets, but tensor-representation models have been recently breaking through. Multi-dimensional arrays, or tensors, are data structures that generalize the concept of vectors and matrices to highly-dimensional domains. In recent years, tensors have also been applied to address numerous data science and machine learning tasks, from simple interpolation to supervised classification \cite{papalexakis2016tensors}.

In this data-science context, a problem of particular interest is that of tensor decomposition, which tries to estimate a set of latent factors that summarize the tensor. Many tensor decompositions were developed as the generalization of well-known matrix-decomposition methods to high-dimensional domains \cite{kolda2009tensor,sidiropoulos2017tensor}. This was the case of the PARAFAC tensor decomposition \cite{harshman1970foundations} and its generalization, the Tucker tensor decomposition \cite{tucker1966some}, which can be both understood as higher-order generalizations of the SVD decomposition of a matrix. More specifically, these decompositions aim at describing (approximating) the tensor as a sum of rank-1 tensors, decomposing it as a sum of outer products of vectors (called factors). The PARAFAC decomposition is conceptually simple and its representation complexity scales gracefully (the number of parameters grows linearly with the rank). The Tucker decomposition enjoys additional degrees of freedom at the cost of greater complexity (exponential dependence of the number of parameters with respect to the rank). Hierarchical tensor decompositions, such as the Tensor Train (TT) decomposition \cite{oseledets2011tensor} or a hierarchical Tucker (hTucker) decomposition \cite{grasedyck2013literature}, try to alleviate this problem. The former unwraps the tensor into a chain of three-dimensional tensors, and the latter generalizes the same idea by organizing the dimensions in a binary tree. Furthermore, in recent years significant effort has been devoted to modify existing decomposition algorithms to deal with factor constraints (e.g., non-negativeness), promote certain priors (e.g., factor sparsity), or be robust to imperfections  \cite{wang2019higher} \cite{xie2017kronecker} \cite{kaya2018parallel}. 

 However, little to no work has been carried out to study the tensor decomposition from a multi-resolution perspective. This can be specially interesting for tensor signals such as videos, where 2-, 3-, and 4-dimensional components are mixed in a single tensor. In this work, we postulate a simple but novel multi-resolution low-rank decomposition method. More specifically, this paper:\vspace{-.15cm}
\begin{itemize}[leftmargin=*]
    \item Introduces a new multi-resolution tensor decomposition to exploit the low-rank structure of a tensor at different resolutions.\vspace{-.2cm}
    \item Proposes an algorithm to implement the decomposition.\vspace{-.2cm}
    \item Tests the benefits of the model via numerical simulations.\vspace{-.15cm}
\end{itemize}

Regarding the first contribution, rather than postulating a low-rank decomposition of the tensor using the original multidimensional representation, we 1) consider a \emph{collection} of lower-order multidimensional representations of the tensor (where several of the original modes of the tensor are combined into a single one); 2) postulate a low-rank decomposition for each of the lower-dimensional representations; 3) map each of the representations back to the original tensor domain; and 4) model the original tensor as the sum of such low-rank representations. As illustrated in detail in the manuscript, this results in an efficient decomposition method capable of combining low-rank structures present at different resolutions. 

Section 2 introduces notation and tensor preliminaries. Section 3 presents our decomposition method. A simple algorithmic approach to address the decomposition is described in Section 4. Illustrative numerical experiments are provided in Section 5. 



\section{Notation and tensor preliminaries} \label{S:Math_Prelim}
\vspace{-.2cm}
The entries of a (column) vector $\bbx$, a matrix $\bbX$ and a tensor $\tenbX$ are denoted by $[\bbx]_n$, $[\bbX]_{n_1,n_2}$ and $[\tenbX]_{n_1,n_2,...,n_I}$, respectively, with $I$ denoting the order of tensor $\tenbX$. Moreover, the $n$th column of matrix $\bbX$ is denoted by $[\bbX]_n$. Sets are represented by calligraphic capital letters. The cardinality of a set $\ccalS$ is denoted by $|\ccalS|$. When a set $\ccalS$ is \emph{ordered}, we use the notation $\ccalS(i)$ with $1\leq i \leq |\ccalS|$ to denote the $i$th element of the set. The vertical concatenation of the columns of matrix $\bbX$ is denoted by $\vect(\bbX)$. $\| \bbX \|_F$  is the Frobenious norm of matrix $\bbX$, which can be equivalently written as $\| \vect(\bbX) \|_2$.

\subsection{Tensor to matrix unfolding}

Given a tensor $\tenbX$ of order $I$ and size $N_1\times ... \times N_I$, there are many ways to unfold the entries of the tensor into a matrix $\bbX$. In this section, we are interested in unfoldings where the columns of matrix $\bbX$ represent one of the original modes of $\tenbX$ and the rows of $\bbX$ represent all the other modes of the tensor. Mathematically, we define the matrix unfolding operator as

\begin{align}\label{E:matrix_tensor_unfolding}
\bbX=&\mathrm{mat}_{p}(\tenbX)\in\reals^{(N_1...N_{p-1}N_{p+1}...N_I)\times N_p}\;\mathrm{where}\hfill\\
 &[\bbX]_{k,n_p}=[\tenbX]_{n_1,...,n_I} \;\text{and}\; \nonumber\\
 &k = n_1+\sum_{i=2,i\neq p}^{I}(n_i-1)\prod_{j=2,j\neq p}^{i-1}N_j.\nonumber
\end{align}
where $p\leq I$ and, to simplify exposition, we have assumed that $p>1$.

\subsection{Tensor to lower-order tensor unfolding}

Consider a tensor $\tenbX$, of order $I$, and let $\ccalI:=\{1,2,...,I\}$ denote the set containing the indexes of all the modes of $\tenbX$. 

\begin{definition}
The ordered set $\ccalP=\{\ccalP_1,...,\ccalP_P\}$ is a partition of the set $\ccalI$ if it holds that: $ \ccalP_p \neq \emptyset$ for all $p$, $ \ccalP_p \cap \ccalP_{p'} = \emptyset $ for all $p' \neq p$, and $\bigcup_{p=1}^P\ccalP_p=\ccalI$. 
\end{definition} 

We are interested in reshaping the entries of the $I$th order tensor $\tenbX$ of size $N_1\times ... \times N_I$ to generate a lower-order tensor $\check{\tenbX}$, with order $P<I$ and according to a given partition $\ccalP=\{\ccalP_1,...,\ccalP_P\}$ as specified next 
\begin{align}\label{E:reshapeten_P}
&\check{\tenbX}=\mathrm{ten}_\ccalP(\tenbX)\in \reals^{\prod_{j=1}^{|\ccalP_1|} |\ccalP_1(j)|\times ... \times \prod_{j=1}^{|\ccalP_P|} |\ccalP_P(j)|}\hfill\\
&\hspace{.5cm}[\check{\tenbX}]_{k_1,...,k_{|\ccalP|}}=[\tenbX]_{n_1,...,n_I}\;\;\mathrm{and}\nonumber\\
&\hspace{.5cm}k_p=n_{\ccalP_p(1)}\;\; \mathrm{if}\;\;|\ccalP_p|=1\nonumber\\
&\hspace{.5cm}k_p= n_{\ccalP_p(1)}+\sum_{i=2}^{|\ccalP_p|}
(n_{\ccalP_p(i)}-1)\prod_{j=1}^{i-1}N_{\ccalP_p(j)}\;\; \mathrm{if}\;\;|\ccalP_p|>1\nonumber
\end{align}
Note that, according to definition of the $\mathrm{ten}_\ccalP(\cdot)$ operator, the indexes along the $p$th mode of $\check{\tenbX}$ represent tuples $(m_{\ccalP_p(1)},...,m_{\ccalP_p(|\ccalP_p|)})$ of indexes of the original tensor $\tenbX$.

Clearly, if $\ccalP=\{\ccalI\}$, so that $P=|\ccalP|=1$ and $|\ccalP_1|=I$, we have that $\mathrm{ten}_\ccalP(\tenbX)=\vect(\tenbX)$. On the other hand, if $\ccalP=\{\{1\},\{2\},...,\{I\}\}$, so that $P=|\ccalP|=I$ and $|\ccalP_p|=1$ for all $p$, we have that $\mathrm{ten}_\ccalP(\tenbX)=\check{\tenbX}$.

Finally, the inverse operator of \eqref{E:reshapeten_P}, which recovers the original tensor $\tenbX$ using as input the reshaped $\check{\tenbX}=\mathrm{ten}_\ccalP(\tenbX)$, is denoted by $\mathrm{unten}_\ccalP(\check{\tenbX})=\tenbX$. Since the definition of $\mathrm{unten}_\ccalP(\cdot)$ starting from \eqref{E:reshapeten_P} is straightforward, it is omitted for conciseness.

\subsection{Low-rank PARAFAC tensor decomposition}
Consider the $I$th order tensor $\tenbX$ along with the matrices $\bbF_i\in \reals^{N_i\times R}$ for $i=1,...,I$. Then, $\tenbX$ is said to have rank $R$ if it can be written as
\begin{equation}\label{E:low-rank-tensor_outerprod_factors}
    \tenbX = \sum_{r=1}^R [\bbF_1]_r\circledcirc [\bbF_2]_r\circledcirc ... \circledcirc   [\bbF_I]_r
\end{equation}
where $\circledcirc$ is the generalization of the outer product for more than two vectors. That is, if $\bbx\in\reals^{N_1}$, $\bby\in\reals^{N_2}$, $\bbz\in\reals^{N_3}$ are three generic vectors, then  $\bbx\circledcirc\bby \circledcirc\bbz $ is a tensor of order $I=3$ satisfying $[\bbx\circledcirc\bby \circledcirc\bbz ]_{n_1,n_2,n_3}=[\bbx]_{n_1}[\bby]_{n_2}[\bbz]_{n_3}\in\reals$. 

The decomposition in \eqref{E:low-rank-tensor_outerprod_factors} is oftentimes referred to as canonical polyadic decomposition or PARAFAC decomposition, with matrices $\bbF_i$ being referred to as factors. As in the case of matrices, moderate values of $R$ induce a parsimonious description of the tensor, since the $\prod_{i=1}^I N_i$ values in $\tenbX$ can be equivalently represented by the $\sum_{i=1}^I RN_i$ entries in $\{\bbF_i\}_{i=1}^I$.  

Using the Khatri-Rao product, denoted as $\odot$, and the different unfolding operators introduced in the previous sections, we have that
\begin{eqnarray}\label{E:low-rank-tensor_matrixized}
\mathrm{mat}_{p}(\tenbX)\!\!\!\! &=&\!\!\!\!\!\! \sum_{r=1}^R \mathrm{mat}_{p}([\bbF_1]_r\circledcirc [\bbF_2]_r\circledcirc ... \circledcirc   [\bbF_I]_r)\nonumber\\
\!\!\!\! &=&\!\!\!\! \!\!(\bbF_I\odot...\odot\bbF_{p+1}\odot\bbF_{p-1}\odot ... \odot\bbF_1) (\bbF_p)^T~~~~
\end{eqnarray}
\begin{eqnarray}\label{E:low-rank-tensor_tensorized}
\mathrm{ten}_{\ccalP}(\tenbX)\!\!&=&\!\!\!\sum_{r=1}^R \mathrm{ten}_{\ccalP}([\bbF_1]_r\circledcirc [\bbF_2]_r\circledcirc ... \circledcirc   [\bbF_I]_r)\nonumber\\
\!\!&=&\!\!\!\sum_{r=1}^R ([\check{\bbF}_1]_r\circledcirc [\check{\bbF}_2]_r\circledcirc ... \circledcirc   [\check{\bbF}_P]_r)\nonumber\\
\mathrm{with}\!\! &~&\! \!\!\check{\bbF}_p=\bbF_{\ccalP_p(|\ccalP_p|)} \odot ... \odot \bbF_{\ccalP_p(2)} \odot \bbF_{\ccalP_p(1)}.~~~~\hfill
\end{eqnarray}
These expressions will be leveraged in the next section.



\section{Multi-resolution low-rank decomposition}
\label{S:mrlr}

Consider a collection of partitions $\ccalP^{(1)}$,...,$\ccalP^{(L)}$, with $|\ccalP^{(l)}|\leq |\ccalP^{(l')}|$ for $l<l'$. Given the $I$th order tensor $\tenbX$ and the collection of partitions $\ccalP^{(1)}$,...,$\ccalP^{(L)}$, we propose the following decomposition for the tensor at hand 
\begin{align}\label{E:multires_low_rank_tensor_decomposition}
   \tenbX &= \sum_{l=1}^L \tenbZ_l, \;\;\mathrm{with}\;\;\mathrm{rank}\big(\mathrm{ten}_{\ccalP^{(l)}}(\tenbZ_l)\big)\leq R_l, \hfill 
\end{align}
which can be equivalently written as
\begin{align}\label{E:multires_low_rank_tensor_decomposition_b}
   \tenbX &= \sum_{l=1}^L \mathrm{unten}_{\ccalP^{(l)}}(\check{\tenbZ}_l),\;\;\mathrm{with}\;\;\mathrm{rank}\big(\check{\tenbZ}_l\big)\leq R_l. \hfill 
\end{align}
where $R_l$ is the rank of the tensor associated to the $l$ partition.

\vspace{.1cm}
\noindent \textbf{Number of parameters:} As already explained, one of the most meaningful implications of low-rank tensor models is the fact that they provide a parsimonious description of the tensor, reducing its implicit number of degrees of freedom. The same is true for the decomposition in \eqref{E:multires_low_rank_tensor_decomposition}. To be concrete, the tensor $\check{\tenbZ}_l=\mathrm{ten}_{\ccalP^{(l)}}(\tenbZ_l)$ has order $P^{(l)}=|\ccalP^{(l)}|$, with the dimension of the $p$th mode being $\prod_{j=1}^{|\ccalP_p^{(l)}|} |\ccalP_p^{(l)}(j)|$. As a result, $\check{\tenbZ}_l$ having rank $R_l$ implies that 
$$R_l \sum_{p=1}^{P^{(l)}} \prod_{j=1}^{|\ccalP_p^{(l)}|} |\ccalP_p^{(l)}(j)| $$  
parameters suffice to fully describe the $\prod_{i=1}^IN_i$ entries in $\check{\tenbZ}_l$.  
Summing across the different $L$ factors implies that
$$\sum_{l=1}^L R_l \sum_{p=1}^{P^{(l)}} \prod_{j=1}^{|\ccalP_p^{(l)}|} |\ccalP_p^{(l)}(j)| $$  
parameters suffice to fully describe the $\prod_{i=1}^IN_i$ entries in $\tenbZ$.



\section{Algorithmic implementation}
\label{S:Algorithmic}

The decomposition introduced in \eqref{E:multires_low_rank_tensor_decomposition} can be obtained by solving the following minimization problem:
\begin{align}
   &\min_{\footnotesize \tenbZ_1 ... \tenbZ_L} \Big\|\tenbX - \sum_{l=1}^L \tenbZ_l\Big\|_F    \label{eq:multiresolution_algorithm}\\
   &\text{s. t.} \;\;\mathrm{rank}\big(\mathrm{ten}_{\ccalP^{(l)}}(\tenbZ_l)\big)\leq R_l. \notag
\end{align}
The approach proposed in this section is to estimate each of the $L$ tensors sequentially, so that when optimizing with respect to $\tenbZ_i$ the remaining tensors $\tenbZ_{l}$ with $l\neq i$ are kept fixed. As a result, the minimization problem to be solved in the $i$th step is:
\begin{align} 
    &\min_{\footnotesize \tenbZ_i} \Big\|\tenbX - \sum_{l\neq i}^L \tenbZ_l - \tenbZ_i\Big\|_F \label{eq:multiresolution_algorithm_b}\\
   &\text{s. t.} \;\;\mathrm{rank}\big(\mathrm{ten}_{\ccalP^{(i)}}(\tenbZ_i)\big)\leq R_i \nonumber
\end{align}
for $i=1,...,L$. The constraint in \eqref{eq:multiresolution_algorithm_b} can be handled using a PARAFAC decomposition %
\begin{equation}\label{eq:ith_tensor_as_parafac}
\tenbZ_i=\sum_{j=1}^{R_i} [\bbH_1^i]_j \circledcirc ... \circledcirc [\bbH_{J_i}^i]_j,
\end{equation}
so that \eqref{eq:multiresolution_algorithm_b} can be equivalently formulated as:
\begin{equation}
    \min_{\bbH_1^i,...,\bbH_{J_i}^i} \Big\|\tenbX - \sum_{l\neq i}^L \tenbZ_l - \sum_{j=1}^{R_i} [\bbH_1^i]_j \circledcirc ... \circledcirc [\bbH_{J_i}^i]_j\Big\|_F.
   \label{eq::multiresolution_algorithm_parafac}
\end{equation}
The above problem is non-convex, but fixing all but one of the factors (say the $j$th one), it becomes linear in $\bbH_j^i$. Under this approach and unfolding the tensor into a matrix $\hat{\bbX}^i=\mathrm{mat}_{p}(\tenbX - \sum_{l\neq i}^L \tenbZ_l)$, we have the following update rule to constructing an Alternating Least Squares (ALS) algorithm:
\begin{equation}
    \min_{\bbH_j^i} ||\hat{\bbX}^i - (\bbH_{J_i}^i \odot ... \odot \bbH_{j+1}^i \odot \bbH_{j-1}^i\odot ... \odot \bbH_1^i)(\bbH_j^i)^T||_F,
   \label{eq:multiresolution_algorithm_parafac_update}
\end{equation}
for all $j=1,...,{J_i}$. Once the ${J_i}$ factors $\{\bbH_j^i\}_{j=1}^{J_i}$ have been obtained, then a) the $i$th tensor $\tenbZ_i$ is found using \eqref{eq:ith_tensor_as_parafac} and b) the problem in \eqref{eq:multiresolution_algorithm_b} is solved for the next $i$, with $i=1,...,L$. As a result, $\sum_{i=1}^L J_i$ instances of   \eqref{eq:multiresolution_algorithm_parafac_update} need to be run. 
Note that, when solving \eqref{eq:multiresolution_algorithm} via \eqref{eq:multiresolution_algorithm_b}-\eqref{eq:multiresolution_algorithm_parafac_update}, the order matters. The first $\tenbZ_l$ to be estimated provides the main (coarser) approximation, while the subsequent ones try to fit the residual error between the main tensor $\tenbX$ and the sum of the previously estimated components $\tenbZ_l$, providing a finer approximation. Due to the structure $\ccalP^{(l)}$, which carries over $\tenbZ_l$, the order in which the tensors $\{\tenbZ_l\}_{l=1}^L$  are approximated is expected to generate variations in the results. 



\subsection{Constructing the partitions}
The algorithm in the previous section assumes that the partitions $\ccalP^{(1)}$,...,$\ccalP^{(L)}$ are given. 
A simple generic approach to design $\ccalP^{(1)}$,...,$\ccalP^{(L)}$ is to rely on a \emph{regular} multiresolution construction that splits the index set $\ccalI=\{1,2,...,I\}$ into smaller sets with the \emph{same cardinality}. More specifically, one can implement a sequential design with $L=I-1$ steps for which, at step $l\in \{1,...,L\}$ we split $\ccalI$ into $l+1$ index sets with (approximately) the same number of elements. The collection of $L=I-1$ partitions $\ccalP^{(1)}$,...,$\ccalP^{(L)}$ is then naturally given by grouping together the sets obtained in each of those steps.  To be more clear, let $\lfloor \cdot \rfloor$ and $\lceil \cdot \rceil$ be the floor and ceil operators and consider the collection of partitions $\ccalP^{(1)}$,...,$\ccalP^{(L)}$ with $L=I-1$ and where the $l$th element is given by  
\begin{eqnarray}\label{E:Sequential_Partition_Partition_l}
\ccalP^{(l)}\!\!\!\!&=&\!\!\!\!\big\{\ccalP^{(l)}_n\big\}_{n=1}^{l+1},\;\text{with}\\
&~&\ccalP^{(l)}_n\!=\!\big\{\lceil (n-1) I/(l+1) \rceil,..., \lfloor n I/(l+1) \rfloor \big\}\nonumber.
\end{eqnarray}
In the above definition we have adopted the convention that, if $x$ is a whole positive number, $\lfloor x \rfloor = x$ and $\lceil x \rceil = x+1$.  
Clearly, the partition design in \eqref{E:Sequential_Partition_Partition_l} is regular in the sense that it achieves $|\ccalP^{(l)}|=l+1$ for all $l$ and $|\ccalP^{(l)}(n)|\approx I/(l+1)$ for $n=1,...,l+1$. 

To gain insights, suppose for simplicity that our tensor $\tenbX$ of order $I$ has size $\eta \times ... \times \eta$, i.e., that the value of $N_i$ is the same across modes, then the number of parameters required to represent $\tenbX$ using the model in \eqref{E:multires_low_rank_tensor_decomposition} and the partitions in \eqref{E:Sequential_Partition_Partition_l} is approximately
\begin{eqnarray}\label{E:NumbeParams_Sequential_Partition}
\sum_{l=1}^{I-1} R_l (l+1) \eta ^{I/(l+l)},
\end{eqnarray}
which contrasts with the $\prod_{i=1}^I N_i=\eta^I$ entries in $\tenbX$.

Clearly, alternative ways to build the partitions $\ccalP^{(1)}$,...,$\ccalP^{(L)}$ are possible. This is especially relevant when prior knowledge exists and one can leverage it to group indexes based on known (di-)similarities among the underlying dimensions. Due to space limitations discussing such alternative partition techniques is out of the scope of this manuscript, but it is part of our ongoing work. 



\section{Numerical experiments}
\label{S:simulations}

The multi-resolution low-rank (MRLR) tensor decomposition scheme is numerically tested in three different scenarios: the first dealing with an amino acids dataset \cite{bro1998multi}, the second one with a video signal \cite{Rozada2021}, and the third one to approximate a multivariate function. The amino acids dataset is a three-mode tensor of size $5 \times 201 \times 61$. The video signal is composed of 173 frames of $1080 \times 720$ pixels each and three channels (R, G, and B). To reduce the computational and memory complexity requirements of the problem, the frames have been sub-sampled and the resolution has been lowered, resulting in a final four-mode tensor of size $9 \times 36 \times 54 \times 3$.  Finally, the multidimensional function in the last scenario has $\reals^3$ as its domain, with each of the three dimensions being discretized using 100 points, so that a tensor with $10^6$ entries is obtained. The Tensorly Python package is used to benchmark the MRLR tensor decomposition against other tensor decomposition algorithms \cite{kossaifi2016tensorly}.

The amino acids tensor $\tenbX$ is approximated using a hierarchical structure of a matrix plus a three-mode tensor. The matrix can be build by unfolding the $5 \times 201 \times 61$ tensor in different ways. Here, two reshapes have been studied, a $201 \times 305$ unfolding (res-1), and a $1005 \times 61$ unfolding (aka res-2). The structure of the algorithm resembles that of a gradient-boosting-like approach \cite{friedman2001greedy}. First, the initial tensor is approximated by a low-rank structure. Then, the residual is approximated by a low-rank structure too. Subsequent residuals are also approximated if necessary. This sequential process can be started from the coarser unfolding, the matrix, or the other way around (reverse). In this experiment, both alternatives have been tested. The rank of the matrix unfolding is fixed while the rank of the three-mode tensor is gradually increased.

\begin{figure}
    \centering
    \includegraphics[width=0.885\linewidth]{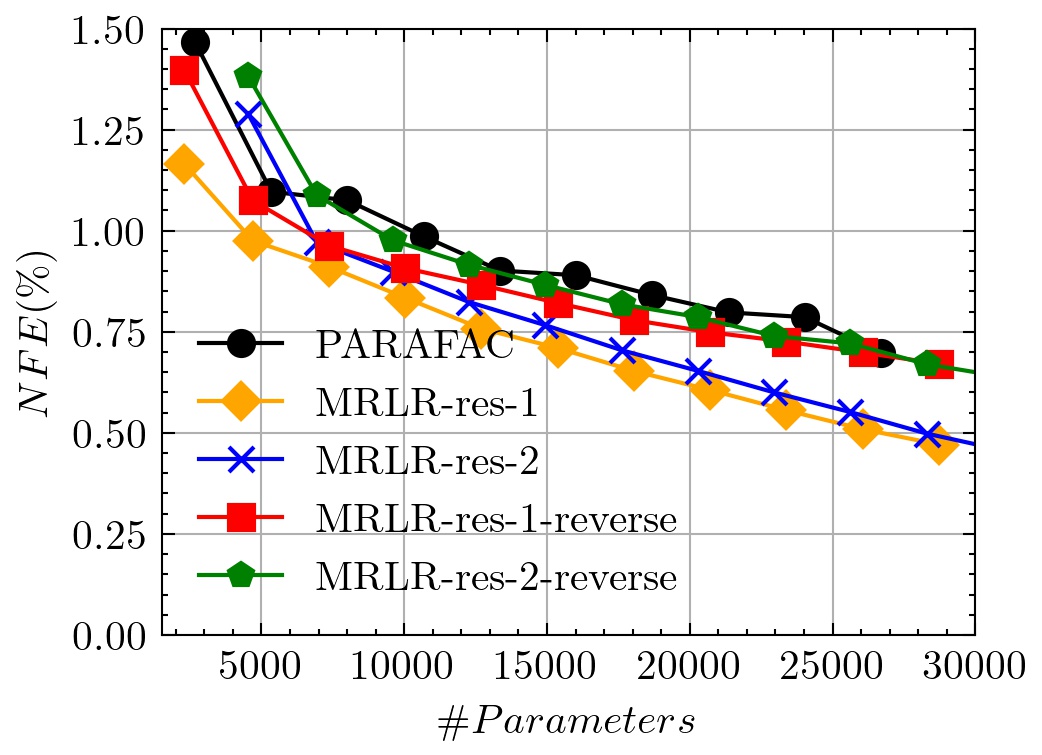}
    \vspace{-.35cm}
    \caption{Normalized Squared Frobenius Error \eqref{eq:NFE} between the original $5 \times 201 \times 61$ amino acids tensor and its approximation obtained via the MRLR and the PARAFAC tensor decompositions when the number of parameters (tensor rank) is increased.}
    \label{fig:fig1}
\end{figure}

\noindent The performance of the algorithms has been measured in terms of Normalized Frobenius Error ($\mathrm{NFE}$) between the true tensor $\tenbX$ and the approximation $\check{\tenbX}$, which is given by
\begin{equation}
    \mathrm{NFE}={||\tenbX - \check{\tenbX}||_F}\big/{||\tenbX||_F}.
    \label{eq:NFE}
\end{equation}

\noindent The results are reported in Fig. \ref{fig:fig1}. The MRLR decomposition is compared to the PARAFAC  decomposition. The res-1 unfolding of the matrix (square-like unfolding) seems to perform better than the res-2 unfolding (tall unfolding). Then, the approximation from the coarser to the finer arrangement beats the reverse one. Moreover, all the MRLR schemes outperform the PARAFAC one in terms of $\mathrm{NFE}$ for the same number of parameters. Indeed, the best-performing MRLR algorithm obtains roughly the same $\mathrm{NFE}$ as the PARAFAC decomposition using $10,000$ parameters less approximately.

\begin{figure}
    \centering
    \includegraphics[width=0.8\linewidth]{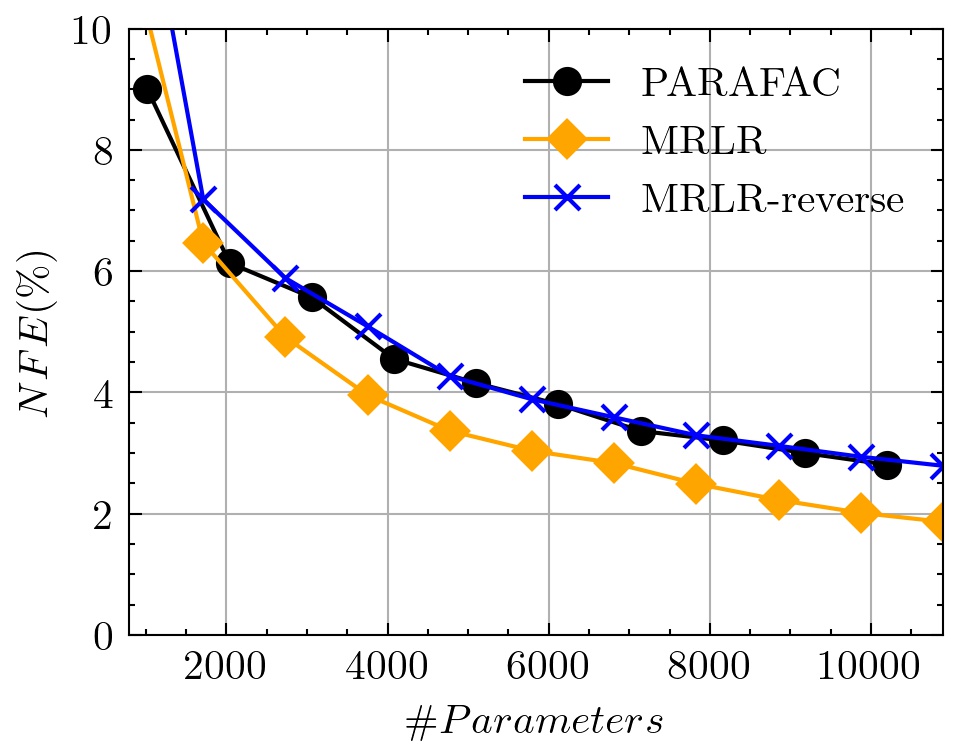}
    \vspace{-.35cm}
    \caption{Normalized Squared Frobenius Error \eqref{eq:NFE} between the original $9 \times 36 \times 54 \times 3$ video signal tensor and its approximation obtained via the MRLR and the PARAFAC tensor decompositions when the number of parameters (tensor rank) is increased.}
    \label{fig:fig2}
\end{figure}

In the second test case, the four-mode video tensor $\tenbX$ is unfolded into a $324 \times 162$ matrix and a $9 \times 36 \times 162$ three-mode tensor. The ranks of the matrix and the three-mode tensors have been fixed to 1. The rank of the four-mode tensor approximation is gradually increased. The results are provided in Fig. \ref{fig:fig2}. Again, the coarser-to-finer arrangement outperforms both, the reverse (finer-to-coarser) arrangement, and the PARAFAC decomposition. It needs approximately $1,500$ parameters less to achieve the same $\mathrm{NFE}$.

Finally, we tested the MRLR tensor decomposition in a third test case to approximate a multivariate function. Given a set of $I$ input variables, with $x_i$ denoting the $i$th input variable and $\ccalX_i$ the set of all possible values of $x_i$, we are interested in functions that map any element $(x_1,...,x_I)\in \ccalX$ into a real value. When these functions are discrete, tensors can be used to model them efficiently. Continuous functions can be discretized/quantized. Tensor decomposition methods can then be leveraged for applications such as approximation, or denoising \cite{kargas2021supervised}. In such a context, we tested the MRLR tensor decomposition algorithm to model the following multivariate continuous function $f: \reals^3 \mapsto \reals$:

\begin{equation}
    f(x_1, x_2, x_3) = \frac{x_1^2 + x_2^2}{e^{|x_2 + x_3|}}.
\label{eq:func_approx}
\end{equation}

\begin{figure}
    \centering
    \includegraphics[width=0.84\linewidth]{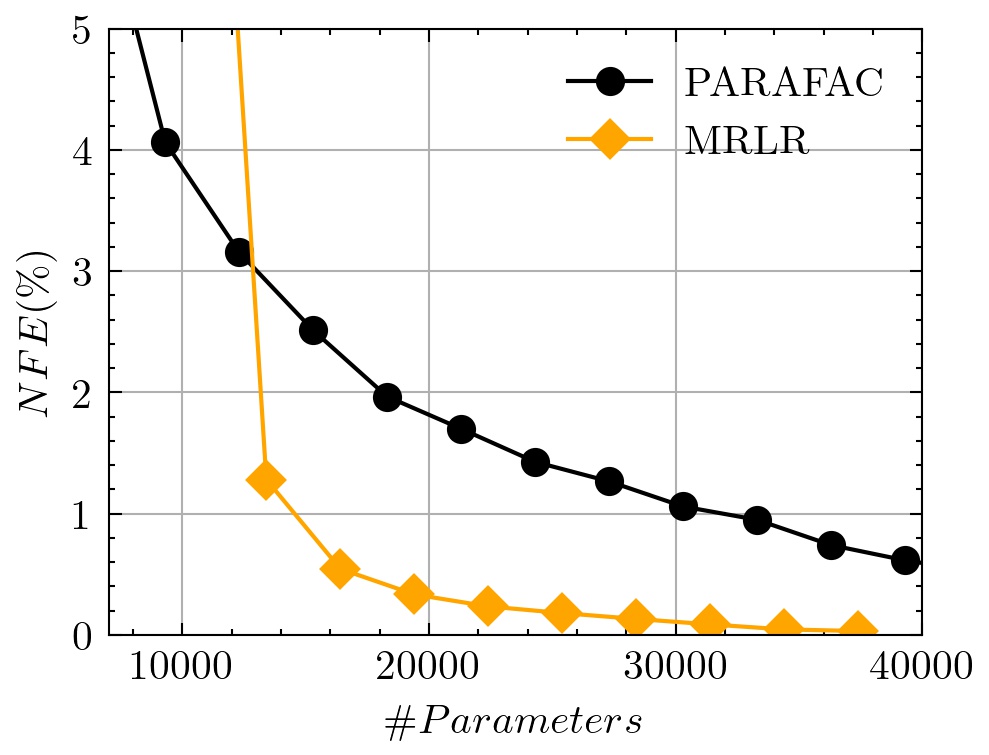}
    \vspace{-.35cm}
    \caption{Normalized Squared Frobenius Error \eqref{eq:NFE} between the $100 \times 100 \times 100$ tensor sampled from the multivariate function in \eqref{eq:func_approx} and its approximation obtained via the MRLR and the PARAFAC tensor decompositions when the number of parameters (tensor rank) is increased.}
    \label{fig:fig3}
\end{figure}

Sampling a three dimensional grid of discrete values ranging from $-5$ to $5$ with an step-size of $0.1$ leads to a $100 \times 100 \times 100$ tensor $\tenbX$ that summarizes the multivariate function in \eqref{eq:func_approx}. The tensor $\tenbX$ can be approximated using the MRLR tensor decomposition to leverage parsimony. The tensor $\tenbX$ is unfolded into a $10000 \times 100$ matrix, and the coarser-to-finer setup has been implemented. The performance of the MRLR tensor decomposition is again compared to that of the PARAFAC decomposition in terms of $\mathrm{NFE}$ for an increasing number of parameters. The results are shown in Fig. \ref{fig:fig3}. As in previous scenarios, the MRLR decomposition outperforms the PARAFAC decomposition for the same number of parameters consistently. At some points, the difference between both algorithms is particularly high. For example, the MRLR tensor decomposition needs roughly $15,000$ parameters to achieve $1\%$ of $\mathrm{NFE}$, while the PARAFAC decomposition needs more than $30,000$ parameters.

\section{Conclusions}
\label{S:conclusions}

This paper presented a parsimonious multi-resolution low-rank (MRLR) tensor decomposition to approximate a tensor as a sum of low-order tensor unfoldings. An Alternating Least Squares (ALS) algorithm was proposed to implement the MRLR tensor decomposition. Then, the MRLR tensor decomposition was compared against the PARAFAC decomposition in two real-case scenarios, and also in a multivariate function approximation problem. The MRLR tensor decomposition outperformed the PARAFAC decomposition for the same number of parameters, showing that it can efficiently leverage information defined at different dimensional orders.


\bibliographystyle{IEEEtran}
\bibliography{references}

\end{document}